\documentclass{timgad}

\usepackage{tikz}
\usepackage{xspace}
\usetikzlibrary{backgrounds,plothandlers,plotmarks,calc,decorations.pathmorphing}

\newcommand{\NN}{\ensuremath{\mathbb{N}}}

\DeclareMathOperator{\im}{im}
\newcommand{\const}{\mathrm{const}}
\newcommand{\KG}{\mathrm{KG}}

\newcommand{\theTitle}{Smaller counterexamples to Hedetniemi's conjecture}
\title{\theTitle}
\author{Marcin Wrochna}
\date{}

\usepackage[backend=bibtex,bibencoding=ascii,style=alphabetic,bibstyle=alphabetic,giveninits=true,doi=false,isbn=false,url=false,maxbibnames=99]{biblatex}
\renewbibmacro{in:}{}
\AtEveryBibitem{\clearname{editor}}

\newbibmacro{string+link}[1]{%
	\iffieldundef{ee}
	{\iffieldundef{url}
		{\iffieldundef{doi}
			{\iffieldundef{eprint}
				{#1}
				{\href{http://arxiv.org/abs/\thefield{eprint}}{#1}}}
			{\href{http://dx.doi.org/\thefield{doi}}{#1}}}
		{\href{\thefield{url}}{#1}}}
	{\href{\thefield{ee}}{#1}}} 
\DeclareFieldFormat{title}{\usebibmacro{string+link}{\mkbibemph{#1}}}
\DeclareFieldFormat[article,inproceedings,inbook,incollection,mastersthesis,thesis]{title}{\usebibmacro{string+link}{\mkbibquote{#1}}}

\begin{document}

\maketitle

\begin{abstract}
	Hedetniemi's conjecture~\cite{hedetniemi1966homomorphisms} for $c$-colorings states that
	the tensor product $G \times H$ is $c$-colorable if and only if $G$ or $H$ is $c$-colorable.
	El-Zahar \& Sauer~\cite{El-ZaharS85} proved it for~$c = 3$.
	In a recent breakthrough, Shitov~\cite{Shitov19} showed counterexamples, for large $c$.
	While Shitov's proof is already remarkably short, Zhu~\cite{Zhu20} simplified the argument and gave a more explicit counterexample for $c=125$.
	Tardif~\cite{Tardif20} showed that a modification of the arguments allows to use ``wide colorings'' to obtain counterexamples for $c=14$, and $c=13$ with a more involved use of lexicographic products.
	This note presents two more small modifications, resulting in counterexamples for $c=5$ (with $G$ and $H$ having 4686 and 30 vertices, respectively).
\end{abstract}

\section{Introduction}
Shitov's counterexample~\cite{Shitov19} to Hedetniemi's conjecture~\cite{hedetniemi1966homomorphisms}
relies on the existence of a graph $F$ which on one hand has high odd girth ($>5$, meaning no cycles of length 3 nor 5), and on the other hand has high fractional chromatic number~($\chi_f>3$).
The second condition means that the chromatic number of the lexicographic product~$F[K_k]$ (``blowing-up'' each vertex into a $k$-clique with all possible edges between adjacent cliques) increases with $k$ as $\chi(F[K_k]) \geq \chi_f(F) \cdot k$ for~all~$k$.
Shitov showed that for any such~$F$, one can find integers $k,c$ and a graph $H$ such that for $G := F[K_k]$ we have $\chi(G)  \geq \lceil \chi_f(F) \cdot k\rceil > c$ and $\chi(H) > c$, yet the product $G \times H$ is $c$-colorable, disproving Hedetniemi's conjecture for $c$-colorings (with $c\simeq 3^{95}$).

As already proved by El-Zahar~\&~Sauer~\cite{El-ZaharS85}, one can assume without loss of generality that $H$ is the exponential graph $K_c^{G}$, or a subgraph thereof, if one seeks to minimise $|V(H)|$ or $|E(H)|$.
Zhu~\cite{Zhu20} simplified Shitov's proof, giving an explicit construction of $H$ (pointing to explicit vertices in the exponential graph, instead of relying on an asymptotic argument on their existence).
He also demonstrated that replacing the second condition directly with $\chi(F[K_k]) > c$ allows to improve parameters,
as $\chi(F[K_k])$ can be strictly larger than the lower bound $\lceil \chi_f(F) \cdot k \rceil$ for small $k$.
This allowed him to give a counterexample $G \times H$ for $c=125$, $|V(G)| = 3403$ and $|V(H)|=10501$.

Tardif~\cite{Tardif20} realized that the high odd girth condition can also be optimized.
For a graph~$G$ and an odd integer $d$, let $\Gamma_d G$ denote the $d$-th power of $G$ (i.e. of its adjacency matrix):
it~has the same vertex set and an edge $uv$ whenever $u$ and $v$ are connected by a walk of length exactly $d$ in~$G$.
An $n$-coloring of $\Gamma_d G$ is also known as a ``$\frac{d+1}{2}$-wide coloring'' of $G$~\cite{SimonyiT06}.
In general, saying $F$ has odd girth $>d$ is equivalent to saying $\Gamma_d F$ has no loop (for odd $d$),
which is equivalent to saying that $\Gamma_d F$ has finite chromatic number:
the identity function is a $|V(F)|$-coloring of $\Gamma_d F$.
It turns out that one can use graphs $F$ with an $n$-coloring of $\Gamma_d F$
where $n$ is much smaller than $|V(F)|$, saving much on the number of colors used in a counterexample.
More precisely, Tardif's proof gives the following interplay of parameters $k,c,n$:

\begin{theorem}[{\cite{Tardif20}}]\label{thm:Tardif}
	Let $k,c,n$ be integers such that $c \geq n + k + 1$ and $c \geq 3k+2$.\\
	Let $F$ be a graph such that $\chi(\Gamma_5 F) \leq n$.
	Then $\chi(K_c^{F[K_k]}) > c$.
\end{theorem}

Hence to give a counterexample it suffices to find $F$ such that $\chi(\Gamma_5 F) \leq n$, yet $\chi(F[K_k]) > c$.
It is is useful to rephrase this first condition in terms of graph homomorphisms (see Section~\ref{sec:prelims} for definitions): let us write $G \to H$ if $G$ admits a homomorphism to $H$.
We have $\chi(\Gamma_d G) \leq n$ if and only if $\Gamma_d G \to K_n$.
It turns out $\Gamma_d$ has a ``right adjoint'' construction $\Omega_d$:
for every graph $H$ there is a graph $\Omega_d H$ such that for all graphs $G$:\vspace*{-.2\baselineskip}
\[\Gamma_d G \to H \text{ if and only if } G \to \Omega_d H.\vspace*{-.2\baselineskip} \]
In particular $\Gamma_d F$ has an $n$-coloring if and only if $F\to \Omega_d K_n$.
Hence to give a counterexample it suffices to find $F$ such that $F \to \Omega_d K_n$, yet $F[K_k] \not\to K_c$ (equivalently, $F \not\to \KG(c,k)$, the Kneser graph).
Thus unless one seeks to minimize the size of $F$, one can always assume $F = \Omega_5 K_n$ to find the best possible parameters $k,c,n$ for which a counterexample using Theorem~\ref{thm:Tardif} exists.

The graph $\Omega_d K_n$ (see Section~\ref{sec:prelims} for explicit constructions) and its chromatic number were first given by Gy{\'{a}}rf{\'{a}}s, Jensen,~\&~Stiebitz~\cite{GyarfasJS04} (for $d=3$, under the name $G_n$) and by Simonyi~\&~Tardos~\cite{SimonyiT06} and independently Baum \& Stiebitz~\cite{BaumS05} (for all odd $d$, under the name $W(\frac{d+1}{2},n)$).
They showed that\vspace*{-.4\baselineskip}
\begin{equation}\label{eq:GJS}
	\chi(\Omega_d K_n) = n\vspace*{-.2\baselineskip}
\end{equation}
(in other words, there exists an $n$-chromatic graph which admits a $d$-wide $n$-coloring, for all $n,d$).
In fact, this follows from the more general statement that the construction $\Omega_d$ (first considered for general graphs by Tardif~\cite{Tardif05} for $d=3$ and by Hajiabolhassan~\cite{Hajiabolhassan09} for all odd $d$)
preserves topological lower bounds on the chromatic number~\cite{Wrochna19}.
(See~\cite{Longueville12,matousek2008using} for an introduction to topological lower bounds).

Therefore the lexicographic product $(\Omega_d K_n)[K_k]$ cannot be $(n+2k-3)$-colored (as otherwise $\Omega_d K_n$ would admit a homomorphism to $\KG(n+2k-3, k)$, which is $(n-1)$-colorable).
Altogether Tardif's proof yields the following:

\begin{corollary}[{\cite{Tardif20}}]
	Let $k,c,n$ be integers such that $c \geq n + k + 1$, $c \geq 3k+2$ and $n+2k-3 \geq c$.\\
	Then $G := (\Omega_5 K_n)[K_k]$ and $H := K_c^G$
	satisfy $\chi(G) > c$, $\chi(H) > c$, yet $\chi(G \times H) \leq c$.
\end{corollary}

Plugging in $k=4$, $c=14$, $n=9$ gives a counterexample $G \times H$ for $14$-colorings,
with $|V(G)|=226980$ and $|V(H)| = 141$ (after taking a suitable subgraph of the exponential graph).

Interestingly, Tardif showed that by replacing $K_k$ with a graph that is $k$-chromatic but not topologically $k$-chromatic (which is only possible for $k \geq 4$) one can construct a counterexample for $13$-colorings, with $|V(G)|=397215$ and $|V(H)|=89$.

Instead, this paper presents other modifications which allow to lower $k$ to 2.
The first is to simply to swap $\Omega_d$ with the lexicographic product.
The second is to repeat a certain argument concerning almost-constant functions in the exponential graph.

\paragraph*{Swapping $\Omega_d$ with the lexicographic product}
The first modification is to replace $(\Omega_d K_n)[K_k]$ with $\Omega_d \left(K_n[K_k]\right) = \Omega_d K_{n \cdot k}$.
In other words, instead of the condition that $\Gamma_d F$ has an $n$-coloring and then using $G := F[K_k]$, there will be a condition that $\Gamma_d G$ is $n \cdot k$-colorable.

This has several advantages.
First, the chromatic number of $\Omega_d \left(K_{n \cdot k}\right)$ is $n \cdot k$,
which is significantly larger than that of $(\Omega_d K_n)[K_k]$ (which is $n+2k-2$ for $k\leq 3$ and $d \geq 5$, according to an unpublished proof by Anna Gujgiczer and G\'{a}bor Simonyi).

Second, whereas the use of $F[K_k]$ required earlier proofs to consider walks of length \emph{at~most}~$d$ in~$F$,
after the modification it suffices to consider walks of length \emph{exactly} $d$.
For a vertex $v$ in a graph, let $N^{\leq d}(v)$, $N^{=d}(v)$ denote the set of vertices reachable from $v$ by a walk of length at most $d$ and exactly $d$, respectively.
The difference is that in a graph of high odd girth, $N^{\leq d}(v)$ is bipartite (2-colorable), while $N^{=d}(v)$ is independent (1-colorable).
An analogous property of wide colorings (see Observation~\ref{obs:local}) 
will allow us to replace $2k$ with $k$ at one place in the proof.

Altogether, this allows to construct counterexamples for $c=7$, as described in Section~\ref{sec:seven}.
Moreover, the modification means that in the final counterexample, one of the graphs can be simply of the form $\Omega_d K_{n'}$, rather than some lexicographic product of it.
Such graphs are relatively well understood, e.g. their chromatic number follows exactly from a topological lower bound and they are similar to other ``spherical'' graphs such as generalized Mycielskians or Schrijver graphs (see for example~\cite{GyarfasJS04,SimonyiT06,BaumS05,Hajiabolhassan09,Wrochna19}).

\paragraph*{Exploiting almost-constant functions more}
The counterexamples in all versions of the proof are constructed by fixing a graph $G$ as discussed above and pointing to vertices of the exponential graph $K_c^G$, which are functions $\colon V(G) \to [c]$.
A significant part is played by almost-constant functions, that is, functions $h$ such that $h(v)$ takes a single value $i$ for almost all $v$ in $G$, except for a small neighborhood.
The proof argues that many such functions must be colored with a minority color, different from $i$.
Section~\ref{sec:five} describes how using more of such functions and repeating essentially the same argument
allows one to reuse color sets that would otherwise need to be disjoint (at the cost of increasing $d$).

Together with the previous modification this allows to give counterexamples for $c=5$.
In other words, $K_5$ is not multiplicative.
This means the only case left open is $c=4$, as far as the multiplicativity of complete graphs $K_c$ is concerned.

\section{Preliminaries}
\label{sec:prelims}
Let $[n] = \{1,\dots,n\}$.
Denote the complete graph on $n$ vertices as $K_n$ and the cycle graph on $n$ vertices as $C_n$.
A \emph{homomorphism} $f \colon G \to H$ is a function $f \colon V(G) \to V(H)$ mapping edges to edges (i.e. $uv \in E(G)$ implies $f(u)f(v) \in E(H)$).
The \emph{odd girth} of a graph is the length of its shortest odd cycle;
equivalently, it is the minimum odd $n$ such that $C_n$ admits a homomorphism to the graph.
The lexicographic product $G[H]$ of two graphs $G,H$ is the graph with vertex set $V(G) \times V(H)$ and with a pair $(g,h)$ adjacent to $(g',h')$ whenever $gg' \in E(G)$ or ($g=g'$ and $hh' \in E(H)$).
The fractional chromatic number $\chi_f(G)$ is the infimum over $\frac{\chi(G[K_k])}{k}$ over $k \in \NN$.

\paragraph*{Wide colorings}
A walk in $G$ is a sequence of vertices where consecutive pairs form by edges of $G$ (not necessarily distinct).
For a graph $G$ and an integer $d$, $\Gamma_d G$ is the graph with the same vertex set $V(G)$ and with an edge between $u,v \in V(G)$ if there is a walk of length exactly $d$ between them in $G$.
In other words, the adjacency matrix of $\Gamma_d G$ is the $d$-th power of the adjacency matrix of $G$.
For a set $S$ of vertices, define $N^{\leq d}(S)$ and $N^{= d}(S)$ to be the set of vertices in $G$ reachable from $S$ by a walk of length at~most~$d$ and exactly~$d$, respectively.
Note that $N^{=d}$ is more than just the set of vertices at distance exactly $d$.
Indeed, $N^{=d}(S)$ is contained in $N^{=d+2}(S)$ for $d \geq 0$ (for $d=0$ this requires that $G$ has no isolated vertices; let us assume this implicitly henceforth).

\begin{observation}\label{obs:local}
	Let $G$ be a graph without isolated vertices, let $\gamma \colon V(G) \to [n]$ be any function and let $d \geq 0$.
	The following are equivalent:

	\begin{enumerate}[label={(\arabic*)}]
	\item $\gamma$ is an $n$-coloring of $\Gamma_{2d+1} G$;
	\item for every color class $S := \gamma^{-1}(i)$ ($i \in [n]$), $N^{=d}(S)$ is independent in $G$;	
	\item for every color class $S := \gamma^{-1}(i)$ ($i \in [n]$) and every $d' \leq d$, $N^{=d'}(S)$ is independent in $G$;
	\item for every color class $S := \gamma^{-1}(i)$ ($i \in [n]$), $N^{\leq d}(S)$ is bipartite in $G$;
	\end{enumerate}
\end{observation}
\begin{proof}
	Suppose $(1)$ holds. Equivalently, for every $u,v$ connected by a walk of length exactly $2d+1$ in $G$, $\gamma(u) \neq \gamma(v)$.
	Let $i \in [n]$ and $S := \gamma^{-1}(i)$.
	If there was any $d' \leq d$ and edge $uv$ in $G$ with $u,v \in N^{=d'}(S)$, then $u,v$ would be connected by walks of length $d'$ to some $i$-colored vertices $u',v'$, respectively; but such $u',v'$ are connected by a walk of length exactly $2d'+1$ and hence by a walk of length exactly $2d+1$ as well, a contradiction.
	Hence $N^{=d'}(S)$ is independent in $G$ for all $d' \leq d$.
	Moreover, $N^{=d}(S)$ is disjoint from $N^{=d-1}(S)$,
	as otherwise for any vertex $u \in N^{=d}(S) \cap N^{=d-1}(S)$ and any neighbor $v$ of $u$ in $N^{=d-1}(S)$ one would have an edge $uv$ in  $N^{=d-1}(S)$.
	Since $N^{\leq d}(S) = N^{= d}(S) \cup N^{=d-1}(S)$, 
	the sets $N^{= d}(S)$ and $N^{=d-1}(S)$ give a bipartition of the subgraph of $G$ induced by $N^{\leq d}(S)$.
	This shows (2),(3),(4).

	Conversely, suppose (4) holds.
	Then $N^{\leq d}(S)$ induces no closed walk of odd length, hence $N^{= d'}(S)$ is independent for $d' \leq d$.
	This implies (3), which trivially implies (2).
	
	Suppose (2) holds.
	Let $u,v$ be any vertices connected by a walk of length exactly $2d+1$ in $G$.
	Let $u'v'$ be the middle edge of such a walk.
	Then it cannot be that $\gamma(u)=\gamma(v)=i$, as otherwise $N^{=d}(\gamma^{-1}(i))$ would contain the edge $u'v'$.
	This shows (1).
\end{proof}

\paragraph*{The $\Omega_{2d+1}$ construction}
For a graph $H$, $\Omega_{2d+1} H$ can be defined up to homomorphic equivalence as the graph such that $\Gamma_{2d+1} G \to H$ if and only if $G \to \Omega_{2d+1} H$ for all $G$.
It can be defined up to isomorphism by requiring that it is a core.
Here are two explicit ways to construct $\Omega_{2d+1} H$~\cite{SimonyiT06,BaumS05,Tardif05}.

In the first construction, for a graph $H$ and an integer $d$,
the vertices of $\Omega_{2d+1} H$ are tuples $(A_0,A_1,\dots,A_d)$ of vertex subsets $A_i \subseteq V(H)$ such that:
\begin{itemize}
	\item $A_0$ is a singleton;
	\item $A_1$ is non-empty;
	\item $A_i \subseteq A_{i+2}$ for $i = 0,\dots,d-2$;	
	\item $A_{d-1}$ and $A_d$ are fully adjacent (that is, $uv \in E(H)$ for all $u \in A_{d-1}$ and $v \in A_d$).
\end{itemize}
The edges are pairs $(A_0,\dots,A_d)$, $(B_0,\dots,B_d)$ such that:
\begin{itemize}
\item $A_i \subseteq B_{i+1}$ and $B_i \subseteq A_{i+1}$ for $i = 0, \dots, d-1$;
\item $A_d$ and $B_d$ are fully adjacent.
\end{itemize}
(If one skips the last three conditions in the definition of a vertex, one gets the same graph plus some isolated vertices, which is of course homomorphically equivalent).
\smallskip

The second construction applies to the case $H = K_n$. While it appears quite different, it is easily shown to be isomorphic to the above.
The vertices of $\Omega_{2d+1} K_n$ are tuples $(x_1,\dots,x_n) \in \{0,\dots,d+1\}^n$ such that:
\begin{itemize}
\item there exists exactly one $i \in [n]$ for which $x_i=0$;
\item there exists some $i \in [n]$ for which $x_i=1$.
\end{itemize}
The edges are pairs $(x_1,\dots,x_n)$, $(y_1,\dots,y_n)$ such that for all $i \in [n]$,
either $|x_i-y_i| = 1$ or $x_i=y_i=d+1$.
From the second construction it is easy to see that $|V(\Omega_{2d+1} K_n)| = n \cdot \left((d+1)^{n-1} - d^{n-1}\right)$.


\paragraph*{Exponential graphs}
As in earlier proofs, we will describe counterexamples by giving a graph $G$ and pointing to vertices of the exponential graph $K_c^G$.
Recall that the vertices of $K_c^G$ are the functions $V(G) \to [c]$,
and two such functions $f,g$ are adjacent if $f(v) \neq g(v')$ for all $vv' \in E(G)$.
The product $G \times K_c^G$ admits a $c$-coloring, given by $(v,f) \mapsto f(v)$.
In fact $K_c^G$ is the most general such graph: for any graph $H$, we have $G \times H \to K_c$ if and only if $H \to K_c^G$.
Thus to give a counterexample for $c$-colorings,
it suffices~\cite{El-ZaharS85} to look for $G$ such that $\chi(G) > c$ and $\chi(K_c^G) > c$.

\section{Counterexample for 7-colorings}\label{sec:seven}

The reader can focus on the case $k = 2$, $c=7$, $n=4$ and $G := \Omega_{5} K_{c+1}$ (so $\chi(G) > c$ by Eq.~\eqref{eq:GJS}).\looseness=-1

\begin{theorem}
	Let $k,c,n$ be integers such that $n \geq k+1$ and $c \geq n+k+1$.
	Let $G$ be a graph such that $\chi(\Gamma_{5} G) \leq n \cdot k$.
	Then $\chi(K_c^G) > c$.
\end{theorem}
\begin{proof}
	Let $\gamma$ be a coloring of $\Gamma_{5} G$ with colors $[n] \times [k]$.
	Let $\alpha,\beta$ be the functions mapping $(a,b) \in [n] \times [k]$ to $a \in [n]$ and $b \in [k]$, respectively.
	By Observation~\ref{obs:local}, $N^{=2}(\gamma^{-1}(a,b))$ is an independent set of $G$ for $(a,b) \in [n] \times [k]$.	
	Suppose there is a $c$-coloring $\Phi$ of $K_c^G$.
	Let us define the following vertices of $K_c^G$, as functions from $V(G)$ to $[c]$:
	\begin{itemize}
	\item $\const_i(v) := i$, for $i \in [c]$, are the constant functions.\\
	A function $f$ is adjacent to $\const_i$ if and only if $i \not\in \im f$.
	Assume without loss of generality $\Phi(\const_i) = i$ (otherwise compose $\Phi$ with a suitable permutation of $[c]$).	
	Then $\Phi(f) \in \im(f)$ for all functions $f$. 
	
	\item $f(v) := \alpha(\gamma(v))$ is a special function that uses all colors in $[n]$.\\
	Hence $\Phi(f) \in [n]$; without loss of generality suppose $\Phi(f) = 1$.\\
	
	\item $h_j(v) := \begin{cases}
			1 & \text{if } v \not\in N^{=1}(\gamma^{-1}(\alpha^{-1}(1))) \\
			j & \text{otherwise}  \end{cases}$
			\qquad for $j \in [c] \setminus [n]$.\\[1ex]
	These are functions that only use two colors: $1$ everywhere except for a small neighborhood.
	Observe that $h_j$ is adjacent to $f$:
	indeed, suppose to the contrary that $h_j(v)=f(v')$ for some $vv' \in E(G)$;
	since $j > n$, it must be that $f(v') = 1$,
	which would imply $v' \in \gamma^{-1}(\alpha^{-1}(1))$,
	hence $v \in N^{=1}(\gamma^{-1}(\alpha^{-1}(1)))$, and thus $h_j(v)=j \neq 1 = f(v')$, a contradiction.
	Therefore $\Phi(h_j) = j$, despite $h_j$ being almost constantly 1.
	
	\item Let $x_1,\dots,x_k$ be distinct colors in $[n]\setminus\{1\}$ (arbitrarily fixed; recall that $n \geq k+1$).\\
	For $j \in [c]\setminus [n]$ let\\[1ex]
	$g_j(v) := \begin{cases}
			j & \text{if } v \not\in N^{=2}(\gamma^{-1}(\alpha^{-1}(1))) \\
			x_b\text{, for some $b \in [k]$ such that }v \in N^{=2}(\gamma^{-1}(1,b)) & \text{if } v \in N^{=2}(\gamma^{-1}(\alpha^{-1}(1)))\\
		\end{cases}$\\[1ex]
	In case there are many such $b$ the choice is arbitrary; one can take the minimum such $b$, say.
	Observe that $g_j$ in the second case gives some color in $\{x_1,\dots,x_k\}$;
	in general $g_j$ only gives colors in $\{x_1,\dots,x_k,j\}$.
	
	We claim that $g_j$ is adjacent to $h_j$, for all $j \in [c] \setminus [n]$.
	Indeed, $1 \not \in \{x_1,\dots,x_k,j\}$,
	so the equality $h_j(v) = g_j(v')$ for an edge $vv'$ of $G$ can only happen if
	$h_j(v) = j$.
	But then $v \in N^{=1}(\gamma^{-1}(\alpha^{-1}(1)))$,
	while $g_j(v') = j$ implies $v' \not\in N^{=2}(\gamma^{-1}(\alpha^{-1}(1)))$,
	so $v$ cannot be adjacent to $v'$, a contradiction.
	Therefore $\Phi(g_j) \neq j$ for all $j \in [c] \setminus [n]$.
	Hence the functions $g_j$ are colored with only $k$ colors only: $\Phi(g_j) \in \{x_1,\dots,x_k\}$.
	
	We now claim that the functions $g_j$ form a clique (of size $c-n$).
	Suppose to the contrary that for some $j \neq j' \in [c] \setminus [n]$ 
	and some $vv' \in G$ we have $g_j(v) = g_{j'}(v')$.
	Since $j$ and $j'$ are distinct and disjoint from $\{x_1,\dots,x_k\}$,
	the common color must be $x_b = g_{i'}(v') = g_i(v)$ for some $b \in [k]$.	
	Then both $v,v' \in N^{=2}(\gamma^{-1}(1,b))$, contradicting the assumption on $\gamma$.
	
	Thus $\Phi$ colors a clique of size $c-n$ with $k$ colors.
	Since we assumed $c -n\geq k+1$, this is a contradiction.\qedhere
	\end{itemize}
\end{proof}

\noindent
Using $G := \Omega_5 K_{c+1}$, so that $\chi(G) > c$ by Eq.~\eqref{eq:GJS}, and assuming $c+1 \leq n \cdot k$ so that $\Gamma_5 G \to K_{c+1} \to K_{n \cdot k}$, we conclude:
\begin{corollary}
	Let $k,c,n$ be integers such that $c \geq n + k + 1$, $n \geq k+1$, and $c + 1 \leq n \cdot k$.
	Then $G := \Omega_5 K_{c+1}$ and $H := K_c^{G}$ satisfy $\chi(G) > c$, $\chi(H) > c$, yet $\chi(G \times H) \leq c$.
\end{corollary}

Plugging in $k = 2$, $c=7$, $n=4$ one gets a counterexample for $7$-colorings (in other words, $K_7$ is not multiplicative).
Inspecting the proof, one can observe that as the graph $H$, it suffices to take the subgraph of $K_c^{G}$ consisting of vertices $\const_i$ ($i \in [c]$), $f$, and then for each possible value of $q := \Phi(f) \in [n]$,
vertices $h^q_j$ and $g^q_j$ ($j \in [c] \setminus [n]$, defined by replacing $1$ with $q$).
Thus $G$ in the counterexample is a graph with $(c+1) \cdot (3^c - 2^c) = 16472$ vertices,
while $H$ has $c + 1 + n \cdot (c-n) \cdot 2 = 32$ vertices.

\pagebreak

\section{Counterexample for 5-colorings}\label{sec:five}
Below one can use $k = 2$, $c=5$, $n=3$ and $G := \Omega_{13} K_{6}$ (so $\chi(\Gamma_{13} G) \leq n \cdot k$ yet $\chi(G) > c$).

\begin{theorem}
	Let $k,c,n$ be integers such that $c \geq n+1$, $c \geq 2k+1$, and $c \geq 5$.
	Let $G$ be a graph such that $\chi(\Gamma_{13} G) \leq n \cdot k$.
	Then $\chi(K_c^G) > c$.
\end{theorem}
\begin{proof}
	Let $\gamma$ be a coloring of $\Gamma_{13} G$ with colors $[n] \times [k]$.
	Let $\alpha,\beta$ be the functions mapping $(a,b) \in [n] \times [k]$ to $a \in [n]$ and $b \in [k]$, respectively.
	By Observation~\ref{obs:local}, $N^{=6}(\gamma^{-1}(a,b))$ is an independent set of $G$ for $(a,b) \in [n] \times [k]$.	
	Suppose there is a $c$-coloring $\Phi$ of $K_c^G$.
	Let us define the following vertices of $K_c^G$, as functions from $V(G)$ to $[c]$:
	\begin{itemize}
	\item $\const_i(v) := i$, for $i \in [c]$, are the constant functions.\\
	A function $f$ is adjacent to $\const_i$ if and only if $i \not\in \im f$.
	Without loss of generality $\Phi(\const_i) = i$.
	Hence $\Phi(f) \in \im(f)$ for all functions $f$. 
	
	\item $f(v) := \alpha(\gamma(v))$ is a special function that uses all colors in $[n]$.\\
	Hence $\Phi(f) \in [n]$; without loss of generality suppose $\Phi(f) = 1$.
	
	\item $h^{(d)}_{i,j}(v) := \begin{cases}
			i & \text{if } v \not\in N^{= d}(\gamma^{-1}(\alpha^{-1}(1))) \\
			j & \text{otherwise}  \end{cases}$
			\qquad for $i,j \in [c]$, $d \in [5]$.
			\\
	These are functions that only use two colors: $i$ everywhere except for a small neighborhood.
	Observe $h^{(1)}_{1,j}$ is adjacent to $f$ for $j > n$, because $f(v) = 1$ only for $v \in \gamma^{-1}(\alpha^{-1}(1))$.
	Hence $\Phi(h^{(1)}_{1,j}) = j$ for $j \in [c]\setminus[n]$, despite $h$ being almost constantly 1.\\[0.6ex]
	
	We claim $h^{(d)}_{i,j}$ is adjacent to $h^{(d+1)}_{i',j'}$
	assuming $i \neq i'$, $j \neq j'$, $i \neq j'$ (so only $j=i'$ is allowed).
	Indeed, this could only fail if for some edge $vv'$ of $G$ we had
	\[ h^{(d)}_{i,j}(v) = j = i' = h^{(d+1)}_{i',j'}(v'). \]
	But in that case $v \in N^{= d}(\gamma^{-1}(\alpha^{-1}(1)))$ while $v' \not \in N^{= d+1}(\gamma^{-1}(\alpha^{-1}(1)))$, a contradiction.
	We use this to show $\Phi(h^{(d)}_{i,j}) = j$ for various $i,j$ and increasing $d$: at each step the possibilities for $i,j$ get larger.
	
	Thus $h^{(1)}_{1,c}$ is adjacent to $h^{(2)}_{c,i}$ for all $i \in [c] \setminus \{1,c\}$.
	Hence $\Phi(h^{(2)}_{c,i}) = i$ for $i \in [c] \setminus \{1,c\}$.
	
	Thus $h^{(2)}_{c,i}$ is adjacent to $h^{(3)}_{i,j}$ for $i \in [c] \setminus \{1,c\}$, $j \in [c] \setminus \{c, i\}$ and $\Phi(h^{(3)}_{i,j}) = j$ for such $i,j$.
	
	For $j \in [c] \setminus \{c\}$ and $\ell \in [c] \setminus \{j\}$, let $i_{j,\ell}$ be an arbitrary color in $[c] \setminus \{1,j,\ell,c\}$ (such a color exists because $c \geq 5$).
	Then $h^{(3)}_{i_{j,\ell}\,,\,j}$ is adjacent to $h^{(4)}_{j,\ell}$.
	Hence $\Phi(h^{(4)}_{j,\ell}) = \ell$ for all $j \in [c] \setminus\{c\}$ and $\ell \in [c] \setminus \{j\}$.
	
	For $\ell,i \in [c]$ with $i \neq \ell$,
	let $j_{\ell,i}$ be arbitrary in $[c] \setminus \{c,\ell,i\}$. 
	Then $h^{(4)}_{j_{\ell,i}\,,\,\ell}$ is adjacent to $h^{(5)}_{\ell,i}$.
	Hence $\Phi(h^{(5)}_{\ell,i}) = i$ for all $\ell,i \in [c]$ with $\ell \neq i$.	 	
	
	\item For $i \in [c] \setminus [k]$ let
	\[g_i(v) := \begin{cases}
			i & \text{if } v \not\in N^{= 6}(\gamma^{-1}(\alpha^{-1}(1))) \\
			b \in [k] \text{ such that } v \in N^{=6}(\gamma^{-1}(1,b))& \text{if } v \in N^{=6}(\gamma^{-1}(\alpha^{-1}(1)))\\
		\end{cases}\]
	If there are many such $b$ the choice is arbitrary.
	Observe that $g_i$ in the second case gives some color in $[k]$;
	in general $g_i$ only gives colors in $[k] \cup \{i\}$.
	
	We claim that for all $i \in [c] \setminus [k]$, $g_i$ is adjacent to $h^{(5)}_{\ell,i}$ for some $\ell \neq i$ in $[c]$.
	Indeed, one can take an arbitrary $\ell \in [c] \setminus ([k] \cup \{i\})$;
	such $\ell$ exists because $c \geq k+2$.
	Then the equality $h^{(5)}_{\ell,i}(v) = g_i(v')$ for an edge $vv'$ of $G$ can only happen if
	$h^{(5)}_{\ell,i}(v) = i$; but then $v \in N^{=5}(\gamma^{-1}(\alpha^{-1}(1)))$,
	while $g_i(v') = i$ implies $v' \not\in N^{= 6}(\gamma^{-1}(\alpha^{-1}(1)))$,
	so $v$ cannot be adjacent to $v'$, a contradiction.
	Therefore, for all $i \in [c] \setminus [k]$, we have $\Phi(g_i) \neq i$.
	Hence the functions $g_i$ get colored with only $k$ colors.
	
	We now claim that the functions $g_i$ form a clique (of size $c-k$).
	Suppose to the contrary that for some $i \neq i' \in [c] \setminus [k]$ 
	and some $vv' \in G$ we have $g_i(v) = g_{i'}(v')$.
	Then the common color must me $b = g_i(v) = g_{i'}(v')$ for some $b \in [k]$.	
	Then both $v,v' \in N^{=6}(\gamma^{-1}(1,b))$, contradicting the assumption on $\gamma$.

	Thus $\Phi$ colors a clique of size $c-k$ with $k$ colors.
	Since we assumed $c \geq 2k+1$, this is a contradiction.\qedhere
	\end{itemize}
\end{proof}

\begin{corollary}
	Let $c,k,n$ be integers such that $c \geq n + 1$, $c \geq 2k+1$, $c \geq 5$, and $c +1 \leq n \cdot k$.
	Then $K_c$ is not multiplicative.
\end{corollary}

Plugging in $k = 2$, $c=5$, $n=3$ one concludes that $G := \Omega_{13} K_{6}$ and $H := K_c^G$ gives a counterexample for $5$-colorings.

Moreover, by being more careful about which $h^{(d)}_{i,j}$ to use, one can decrease the maximum $d$ and the number of functions needed.
Specifically, for parameters $k = 2$, $c=5$, $n=3$,
let us consider as before all possibilities for $q := \Phi(f) \in [n]$.
Let $q\oplus 1, q \oplus 2$ be the other elements of $[n]$.
Define \[h^{q,(d)}_{i,j}(v) := \begin{cases}
			i & \text{if } v \not\in N^{= d}(\gamma^{-1}(\alpha^{-1}(q))) \\
			j & \text{otherwise}  \end{cases}\]
and \[g^q_i(v) := \begin{cases}
			i & \text{if } v \not\in N^{=3}(\gamma^{-1}(\alpha^{-1}(q))) \\
			q \oplus (b-1) \text{ for } b \in [k] \text{ such that } v \in N^{=3}(\gamma^{-1}(1,b))& \text{if } v \in N^{=3}(\gamma^{-1}(\alpha^{-1}(q)))\\
		\end{cases}.\]			
Then by the same arguments, a counterexample for 5-colorings can be built from 
$G := \Omega_7 K_6$ and the subgraph $H$ given by the following vertices of $K_5^G$:
\begin{itemize}
\item $\const_i$ for $i \in [c]$,
\item $f$; and then for all $q \in [n]$:
\item $h^{q,(1)}_{q,4}$, $h^{q,(1)}_{q,5}$,
\item $h^{q,(2)}_{4,5}$, $h^{q,(2)}_{5,4}$, $h^{q,(2)}_{5,q \oplus 2}$,
\item $g^q_4$, $g^q_5$, $g^q_{q \oplus 2}$.
\end{itemize}
Thus $G$ has $6 \cdot \left(4^{5} - 3^{5}\right) = 4686$ vertices (and 36015 edges),
while $H$ has $c + 1 + n \cdot (2 + 3 + 3) = 30$ vertices (and 108 edges).

\medskip
We remark that for $c=4$ one can also find $d$ large enough such that a similar proof shows $\Phi(h^{(d)}_{i,j}) = j$ for all color pairs $i\neq j \in [c]$ (indeed, the requirement $c\geq 5$ can be removed by using $h^{(4)}_{j,\ell}$ only for $\ell=c$, and then adding a few more similar steps).
So the bottleneck in the proof is the final argument which relies on $c \geq 2k+1$.

%

%

\pagebreak

\printbibliography

\end{document}